# Factional Brownian motion with multivariate time in a large convex area: persistence exponents

## G. Molchan

[1] Institute of Earthquake Prediction Theory and Mathematical Geophysics, Russian Academy of Science, Moscow

[2] National Research University Higher School of Economics, Moscow.

**Abstract**. The fractional Brownian motion of index 0 < H < 1, H-FBM, with d-dimensional time is considered on an expanding set TG, where G is a bounded convex domain that contains 0 at its boundary. The main result: if 0 is a point of smoothness of the boundary, then the log-asymptotics of probability that H-FBM does not exceed a fixed positive level in TG is (H − d + o(1)) log T, T>>1. Some generalizations of this result to isotropic but not self-similar Gaussian processes with stationary increments are also considered.

## 1. Introduction.

This note is the sequel to ( Molchan-1, 2017) where we considered the log-asymptotics of probability

$$p_T(\Delta) = P\{\xi(t) < 1, t \in T\Delta\}, \quad T \to \infty, \quad (1.1)$$

where $\xi(t) = w_H(t), t \in R^d$ is the fractional Brownian motion, FBM, with the Hurst parameter 0<H<1, and $\Delta \subset R^d$ is a ball containing 0 at its boundary. FBM is a centered Gaussian random process with multivariate time and the correlation function

$$Ew_H(t)w_H(s) = 0.5(|t|^{2H} + |s|^{2H} - |t-s|^{2H}) ;$$

$w_H(t)$ is self-similar, isotropic and has stationary increments, i.e. the following

$$\{w_H(\lambda Ot + t_0) - w_H(t_0)\} \doteq \{|\lambda|^H w_H(t)\} \quad (1.2)$$

holds in the sense of the equality of finite-dimensional distributions provided that $\lambda \in R^1$, $O \in SO(d)$, and $t_0 \in R^d$ are fixed.

The result for the one dimensional time and $\Delta = (0,1)$ is well known (Molchan, 1999):

$$\log p_T(\Delta)/\log T = H - 1 + o(1), \quad T \to \infty. \quad (1.3)$$

Aurzada et al (2016) made more precise the remainder term o (1) and extended the result (1.3) on some processes with stationary increments. The papers by Dembo et al (2013) and Aurzada et al (2016) showed the essential role of the time discretization for the persistence probability analysis. The current state of the persistence problem can be found in the surveys (Bray et al, 2013; Aurzada and Simon, 2015).

The persistence exponent (the main term in the asymptotics of type (1.3)) depends on the position where the considered process is equal to 0. According to (Molchan, 1999), if this position is an internal point of $\Delta$, then for a large class of self-similar processes

$$\log p_T(\Delta)/\log T = -d + o(1).$$

In the report (Molchan, 2005), the author put forward the hypothesis that for $\xi(t) = w_H(t), t \in R^d$ and $\Delta = [0,2] \times [-1,1]^{d-1}$

$$\log p_T(\Delta)/\log T = H - d + o(1). \tag{1.4}$$

Partially this hypothesis was confirmed in (Molchan-1, 2017) where (1.4) was proved for any spherical domain $\Delta$ containing 0 at its boundary. Our goal is to prove the following

**Proposition 1**. Suppose that $\Delta$ is a d-dimensional bounded convex domain that contains 0 at its boundary $\partial \Delta$. If $\partial \Delta$ is smooth in a vicinity of 0, then (1.4) holds for $\xi(t) = w_H(t), t \in R^d$.

Moreover

**The lower bound.** Let $\xi(t), t \in R^d$, $\xi(0) = 0$ be a centered Gaussian process such that

a) $\xi(t)$ is isotropic and has stationary increments;

b) $E\xi^2(t) < k|t|^{2H}, |t| < \infty$ and $E\xi^2(t) \sim c^2|t|^{2H}$ as $t \to \infty$, where $0 < H < 1$.

Assume, $\{t \in R^d : t^{(1)} = 0\}$ is the tangent plane to $\Delta$ at $0 \in \Delta$. Then

$$T^{H-d}(\ln T)^{-\lambda} \leq P\{\xi(t) < 1, t \in T\Delta \cap Z^d)$$

with any $\lambda > 1/2$. If, moreover, the process is self-similar, i.e. $\xi(t) = w_H(t)$, then

$$T^{H-d}(\ln T)^{-\lambda} \leq P\{\xi(t) < 1, t \in T\Delta\}$$

with any $\lambda > d/(2H)$.

**The upper bound.** Assume in addition to (a,b)) that $\xi(t)$ has an absolutely continuous power spectrum measure $F(d\lambda) = f(\lambda)d\lambda, \lambda \in R^d$ such that $f(\lambda) > g(|\lambda|) > 0$, where $g(x), x \in R_+^1$ is a bounded continuous decreasing function and

$$\int_0^\infty \frac{|\log g(x)|}{1+x^2} dx < \infty .$$

Then

$$p_T(\Delta) \leq \exp(c\sqrt{\ln T}) T^{-d+H} .$$

Note that for the fraction Brownian motion $E\xi^2(t) = |t|^{2H}$ and $f(\lambda) = c|\lambda|^{-d+2H}$.

The key point of our proof is the following functional

$$F_n = \sum_{i=1}^n (\xi(t_i) - M_{i-1})_+ ,$$

where $t_i$ is a suitable onto-mapping $Z_+^1 \to Z^d$ and $M_n = \max(\xi(t_i), i = 1,...,n)$. A conceptually similar functional was used in (Molchan -2, 2017) to find the persistence exponent for the Integral fractional Brownian motion considered at the domain $T\Delta = (-T, T)$.

## 2. Preliminaries.

Assume that $\Delta$ is a d-dimensional bounded convex domain and 0 is a point of smoothness of its boundary $\partial \Delta$. Then there exist a suitable ball $B_r(t_0)$ and a cube $K_L(s_0)$, such that

$$0 \in B_r(t_0) \subseteq \Delta \subseteq K_L(s_0) . \qquad (2.1$$

Here $r > 0$ and $L$ are the size parameters and $t_0, s_0$ are the centers of the considered figures. In addition, all these three domains have a common tangent plane at 0.
Relation (2.1) implies the following estimates of probability (1.1):

$$p_T(B_r(t_0)) \geq p_T(\Delta) \geq p_T(K_L(s_0)) .$$

Hence to estimate the persistence probability $p_T(\Delta)$ for isotropic processes we only need to consider two central-symmetric domains, $K_1(e) = [-1,1]^d + e$ and $B_1(e)$, with the center $e = (1,0,..,0)$.

In what follows we shall need the following interpolation Lemma. Implicitly, similar statements are often used in persistence probability analysis (see e.g. Aurzada and Monch, 2016; Molchan-1, 2017).

**Lemma 2**. Let $\xi(t), t \in \Delta_T \subset R^d$ be a random process, $U_T = \{t_i\}$ is a $\rho$-net for $\Delta_T$, i.e.

$\cup_i B_\rho(t_i) \supset \Delta_T$; $M_\xi(V) = \sup(\xi(t), t \in V)$ and $|U_T| = \#\{t_i : t_i \in U_T\}$.

If $\xi(t)$ has stationary increments, then for $a > c$ and $b > 0$

$$P(M_\xi(\Delta_T) \leq c) \leq P(M_\xi(U_T) \leq c) \leq P(M_\xi(\Delta_T) \leq a) + |U_T|P(M_\xi(B_\rho(0)) \geq a - c) \quad (2.2)$$

and

$$EM_\xi(U_T) \leq EM_\xi(\Delta_T) \leq EM_\xi(U_T) + b + |U_T|E(M_\xi(B_\rho(0)) - b)_+ \quad (2.3)$$

**Proof**. One has

$$M_\xi(\Delta_T) \leq M_\xi(U_T) + \max_{t_k \in U_T}(M_\xi(B_\rho(t_k)) - \xi(t_k)) := M_\xi(U_T) + \delta_{U_T}.$$

Therefore

$$P(M_\xi(\Delta_T) \leq a) \geq P(M_\xi(U_T) \leq c, \delta_{U_T} \leq a - c) \geq P(M_\xi(U_T) \leq c) - R_{U_T}, \quad (2.4)$$

where

$$R_{U_T} = P(\delta_{U_T} \geq a - c) \leq \sum_k P(M_\xi(B_\rho(t_k)) - \xi(t_k)) \geq a - c).$$

Due to stationary increments of $\xi(t)$, we have

$$R_{U_T} \leq |U_T|P(M_\xi(B_\rho(0)) \geq a - c) . \quad (2.5)$$

Combining (2.4) and (2.5), we get the right-hand part of (2.2). The left-hand part is obvious.

To prove (2.3), note that

$$\delta_{U_T} \leq b + \max_{t_k \in U_T}(M_\xi(B_\rho(t_k)) - \xi(t_k) - b)_+.$$

Hence, the stationarity of $\xi(t)$ increments implies

$$E\delta_{U_T} \leq b + |U_T|E(M_\xi(B_\rho(0)) - b)_+ . \quad (2.6).$$

Now (2.3) follows from (2.6) and from the relation

$$M_\xi(U_T) \leq M_\xi(\Delta_T) \leq M_\xi(U_T) + \delta_{U_T}.$$

Let $\Delta_T = T\Delta, U_T = \Delta_T \cap Z^d := [T\Delta]$ is the set of integer-valued points in $\Delta_T$ .

**Corollary 3.** Let $\xi(t), t \in R^d$, $\xi(0) = 0$ be a centered Gaussian process with stationary increments such that $E\xi^2(t) < c|t|^{2\varepsilon}, |t| < 1$ , for some $\varepsilon > 0$ . Then

$$P(M_\xi([T\Delta]) \leq 0) \leq P(M_\xi(T\Delta) \leq \sqrt{\kappa \ln T}) + CT^{-\nu(\kappa)} , \quad (2.7)$$

$$0 \leq EM_\xi(T\Delta) - EM_\xi([T\Delta]) \leq C\sqrt{\ln T} , \quad (2.8)$$

where $\nu(\kappa) = C_\varepsilon \kappa - d$ .

**Proof.** Relation (2.7) follows from (2.2) with $c = 0$ and $a = \sqrt{\kappa \ln T}$ . Indeed, $|U_T| = |T\Delta| \leq cT^d$ .

Since $E(\xi(t + \delta) - \xi(t))^2 < c|\delta|^{2\varepsilon}, |\delta| < 1$ , the Fernique (1975) inequality for $\xi(t)$ gives

$$P(M_\xi(B_1(0)) \geq rc_\varepsilon) \leq c\int_r^\infty \exp(-x^2/2)dx, \quad r > (1 + 4d)^{1/2}. \quad (2.9)$$

Whence

$$|U_T|P(M_\xi(B_1(0) \geq \sqrt{\kappa \ln T}) \leq cT^{-\kappa/(2c_\varepsilon^2)+d} , \quad T > T_0 ,$$

which proves (2.7) with $\nu(\kappa) = \kappa/(2c_\varepsilon) - d$.

Obviously, (2.9) leads to the moment estimation:

$$|U_T|E(M_\xi(B_1(0)) - \sqrt{\kappa \ln T})_+ \leq cT^{-\nu(\kappa)}, \ T > T_0. \tag{2.10}$$

Now (2.8) follows from (2.3) and (2.10).

**Lemma 4.** Let $\xi(t), t \in R^d$, $\xi(0) = 0$ be a continuous centered Gaussian process with stationary increments such that

$$E\xi^2(t) < k|t|^{2H}, |t| < \infty \tag{2.11}$$

and

$$E\xi^2(t) \sim c^2|t|^{2H} \text{ as } t \to \infty, \tag{2.12}$$

where $0 < H < 1$ and $a \sim b$ means that $a/b \to 1$. Then

$$EM_\xi(T\Delta) \sim cT^H EM_\xi(\Delta), \tag{2.13}$$

where $\Delta$ is a bounded domain containing $0$, $0 \in \Delta$.

**Proof.** The conditions (2.11, 2.12) guarantee the weak convergence of the scaled process $Z_T(t) = T^{-H}\xi(Tt), t \in \Delta$ to the fractional Brownian motion $w_H(t), t \in \Delta$ as $T \to \infty$ on $C([-1,1]^d, R)$ (see Whitt, 2002, Theorems 11.6.5 and 11.6.7).

Due to (2.11), $EZ_T^2(t) \leq k|t|^{2H}$. Applying the Fernique inequality (2.9) to $Z_T(t)$, we get

$$P(M_{Z_T}(\Delta)) \geq rc_H) \leq c\int_r^\infty \exp(-x^2/2)dx.$$

This estimate is uniform in the parameter T. Hence, taking into account the weak convergence of $Z_T(t)$ to $w_H(t)$, we may conclude that

$$\lim_{T\to\infty} EM_{Z_T}(\Delta) = cEM_{w_H}(\Delta).$$

This proves (2.13).

With any centered Gaussian process $\xi(t), t \in \Delta$ is associated the Hilbert space of functions $H_\xi(\Delta)$ with reproducing kernel $K(t,s) = E\xi(t)\xi(s), (t,s) \in \Delta \times \Delta$ ( Lifshits , 2012). By definition the set of functions $\{k_s(t) = K(t,s), s \in \Delta\}$ is dense in $H_\xi(\Delta)$ relative to the norm

$$\left\|\sum a_i k_{s_i}\right\| = \left(\sum K(s_i, s_j)a_i a_j\right)^{1/2}.$$

**Lemma 5.** Let $\xi(t), t \in R^d$, $\xi(0) = 0$ be a centered Gaussian process with stationary increments. Assume that $\xi(t)$ has an absolutely continuous power spectrum measure $F(d\lambda) = f(\lambda)d\lambda, \lambda \in R^d$ such that

$$f(\lambda) > g(|\lambda|) > 0, \tag{2.14}$$

where $g(x), x \in R_+^1$ is a bounded continuous decreasing function, and

$$\int_0^\infty \frac{|\log g(x)|}{1+x^2} dx < \infty \qquad (2.15)$$

Then there exists continuous function $\varphi(t) \in H_\xi(R^d)$ such that $\varphi(0) = 0$ and $\varphi(t) = 1, |t| \geq 1$.

**Proof.** Under conditions (2.14, 2.15) there exists an integer function of the exponential type $U(x), x \in R^1$ such that $U(x) > 0$, $U(x) = U(-x)$ and for any $\varepsilon > 0$

$$\sup_x U(\varepsilon x)/g(x) < \infty. \qquad (2.16)$$

(see Inozemtcev and Marchenko, 1956).

Consider

$$v(\lambda) = U(|\lambda|) - U(|\varepsilon \lambda|), \quad \lambda \in R^d \qquad (2.17)$$

and show that $\psi(\lambda) = v(\lambda)/f(\lambda) \in L_f^2$.

This fact follows from the finiteness of $I_\varepsilon = \int U^2(\varepsilon|\lambda|)/f(\lambda)d\lambda$. One has

$$I_\varepsilon = \int_{|\lambda|\leq 1} [U(\varepsilon|\lambda|)/g(|\lambda|)]^2 [g(|\lambda|)/f(\lambda)]g(|\lambda|)d\lambda +$$

$$\int_{|\lambda|>1} [U\varepsilon|\lambda|)/g(|\lambda|)]^2 [g(|\lambda|)/f(\lambda)]^2 f(|\lambda|)d\lambda.$$

Due to (2.14, 2.16),

$$I_\varepsilon \leq c(\int_{|\lambda|\leq 1} g(|\lambda|)d\lambda + \int_{|\lambda|>1} f(\lambda)d\lambda).$$

The last expression is bounded because $g(x) \leq g(0)$ and $f(\lambda)1_{|\lambda|>\delta} \in L^1(R^d)$ for any $\delta > 0$.

Similarly we conclude that $U(\varepsilon|\lambda|) \in L^1$:

$$\int [U(\varepsilon|\lambda|)d\lambda = \int [U(\varepsilon|\lambda|)/g(|\lambda|)]g(|\lambda|)d\lambda \leq c\int f(\lambda) \wedge g(0)d\lambda < \infty.$$

Because $U(x)$ is an even integer function of the exponential type, the same is true for $U(|\lambda|), \lambda \in R^d$.

Applying the Paley-Wiener theorem (Yosida, 1968), to $U(|\lambda|) \in L^1$ we conclude that $U(|\lambda|)$ is Fourier transform of a continuous finite function $u(t)$.

Due to the spectral representation of the reproducing kernel

$$K(t,s) = E\xi(t)\xi(s) = \int (e^{i(t,\lambda)} - 1)(e^{-i(s,\lambda)} - 1)f(\lambda)d\lambda,$$

the mapping

$$\gamma : H_\xi \ni k_s(t) = K(t,s) \Rightarrow e^{i(t,\lambda)} - 1 \in L_f^2$$

generates the isometric embedding of $H_\xi$ into $L_f^2$. Moreover, the value $k(t)$ of $k(\cdot) \in H_\xi$ is the projection of $\gamma k(\cdot)$ onto $e^{i(t,\lambda)} - 1$ in $L_f^2$.

Hence, projecting $\psi(\lambda) = v(\lambda)/f(\lambda) \in L_f^2$ on the elements $\{e^{i(t,\lambda)} - 1\}$ we get the following function from $H_\xi(R^d)$:

$$\varphi(t) = \int (e^{i(t,\lambda)} - 1)\psi(\lambda)f(\lambda)d\lambda = \int (e^{-i(t,\lambda)} - 1)(U(|\lambda|) - U(|\varepsilon\lambda|))d\lambda$$

$$= (2\pi)^d (u(t) - u(t/\varepsilon)\varepsilon^{-d} + u(0)(\varepsilon^{-d} - 1)).$$

Because $u(t)$ is finite, we have: $\varphi(t) = c$ for $|t| > a$ with some constants $c, a$. Moreover, $c > 0$ for any $\varepsilon < 1$ because $U(x) > 0$, and therefore $(2\pi)^d u(0) = \int U(|\lambda|)d\lambda > 0$.

It is easy to see, that the desired function with parameters $c = a = 1$ is $\varphi(t/a)/c$.

## 3. Proof of the main result.

### The lower bound

As we mention above to estimate $p_T(\Delta) = P\{\xi(t) < 1, t \in T\Delta\}$ from below it is sufficient to consider the domain $\Delta = [-1,1]^d + e$, $e = (1,0,..,0)$.

*Step 1.* Let us consider the process $\xi(t)$ at integer points $t \in Z^d$. We shall need a suitable onto-mapping $t_i: Z_+^1 \to Z^d$ with the following properties:

The curve $i \to t_i$ should be 'continuous':

$$1 \leq |t_i - t_{i+1}| \leq L \tag{3.1}$$

and should consistently cover the surfaces $S_n = \{t : |t| = n\}$ of different levels, i.e. the following is implied:

$$|t_i| < |t_j| \Rightarrow i < j,$$

where $|t| = \max_{\alpha=1 \div d} |t^{(\alpha)}|$.

Let us divide the surface $S_n$ into two disjoint zones $S_n^1$ and $S_n^2$. Consider two points $t_i \in S_n^1$ and $t_j \in S_n^2$ such that $t_i \neq t_k$ for any $k < i$ and $t_j \neq t_m$ for any $m < j$. In this case, we want to have $i < j$. In other words, at first the curve $i \to t_i$ covers the zone $S_n^1$ and then $S_n^2$. The curve may exit from the second zone to the first, but not vice versa.

To specify the zones, note that $S_n$ consists of (d-1)-dimensional faces $\{S_{n\beta}, \beta = 1 \div 2d\}$. If $t_{n\beta}$ is the center of $S_{n\beta}$ and

$$\hat{S}_{n\beta} = \{t \in S_{n\beta} : |t - t_{n\beta}| \leq n(1 - 1/L_n)\} \tag{3.2}$$

where $L_n > 1$ is a slowly varying function, then $S_n^2 = \bigcup_{\beta=1}^{2d} \hat{S}_{n\beta}$ by definition.

If $d > 2$, $S_n^1$ as the complement of $S_n^2$ is a connected zone. Therefore the above-described way of numbering of the lattice vertices, $t \in Z^d$, is quite feasible with given $S_n^1$ and $S_n^2$. In the case $d = 2$, $S_n^1$ consists of 4 disjoint pieces. Therefore our curve has to connect them by means of an exit on a surface of level $n-1$, $S_{n-1}$.

*Step 2.* With the above numbering of $t \in Z^d$ we consider the following functional

$$F_n = \sum_{i=1}^n (\xi_i - M_{i-1})_+ , \tag{3.3}$$

where $\xi_i = \xi(t_i)$, $M_n = \max(\xi_i, i = 1, ..., n)$, $M_0 = 0$.

It is obvious that

$$F_n = M_n - M_0 = M_n, \tag{3.4}$$

$$(\xi_i - M_{i-1})_+ \leq (\xi_i - \xi_{i-1})_+ , \tag{3.5}$$

$$(\xi_i - M_{i-1})_+ = 0 \text{ if } t_i = t_j \text{ and } i > j . \tag{3.6}$$

Hence, from all preimages of a fixed point t, the functional $F_n$ takes into account only one of them, namely, the smallest index associated with that point.

Let $N_n$ be the ordinal number of the last point $t$ of the level $n$. For simplicity we introduce a new symbol:

$$F_{N_n} = F(n\Delta), \quad \Delta = [-1,1]^d,$$

Obviously

$$M_{N_n} = M_\xi([n\Delta]) .$$

As above, $M_\xi(A) = \sup\{\xi(t), t \in A\}$ and $[A] = A \cap Z^d$.

One has

$$F(n\Delta) - F((n-1)\Delta) = \sum_{t_i \in S_n^1} (\xi_i - M_{i-1})_+ + \sum_{t_i \in S_n^2} (\xi_i - M_{i-1})_+ := \sum_n^1 + \sum_n^2$$

It is obvious that

$$\sum_n^2 - \sum_n^1 = M_\xi([n\Delta]) - M_\xi([(n-1)\Delta]) \vee M_\xi(S_n^1) \geq 0 .$$

Therefore

$$F(n\Delta) - F((n-1)\Delta) \leq 2\sum_{t_i \in S_n^2} (\xi_i - M_{i-1})_+ . \tag{3.7}$$

Consider $0 < q < 1$ and $m = [qn]$. By (3.4),

$$F(n\Delta) - F(m\Delta) = M_\xi([n\Delta]) - M_\xi([m\Delta]) .$$

Using (2.8, 2.13), one gets

$$E(F(n\Delta) - F(m\Delta)) \geq C_H n^H - c\sqrt{\ln n} - m^H EM_\xi(\Delta)$$

$$= C_H n^H (1 - q^H) - c\sqrt{\ln n} . \qquad (3.8)$$

On the other hand, due to (3.7), we have

$$E(F(n\Delta) - F(m\Delta)) \leq 2\sum_{k=m}^{n} E\sum_{t_i \in S_k^2} (\xi_i - M_{i-1})_+ .$$

Therefore

$$E(F(n\Delta) - F(m\Delta)) \leq cn^d \max(E(\xi_i - M_{i-1})_+, t_i \in S_k^2, m \leq k \leq n) . \qquad (3.9)$$

Since (3.5), one has

$$E(\xi_i - M_{i-1})_+ \leq E(\xi_i - \xi_{i-1})_+ 1_{\xi_i \geq M_i} \leq \sqrt{\kappa \ln n} P(\xi_i - M_i \geq 0) + R_n , \qquad (3.10)$$

where

$$R_n = E(\xi_i - \xi_{i-1} - \sqrt{\kappa \ln n})_+ .$$

Using (3.1) and stationarity of the increments of $\xi(t)$, we get

$$R_n \leq E(M_\xi(B_L(0)) - \sqrt{\kappa \ln n})_+$$

By (2.10), we can continue

$$n^d R_n \leq cn^{-\nu(\kappa)}, \qquad (3.11)$$

where $\nu(\kappa)$ is a linear function of $\kappa$.

Combining (3.8-3.11) we get

$$cn^{H-d} \leq \sqrt{\ln n} \max\{P(\xi_i - M_i \leq 0), t_i \in S_k^2, m \leq k \leq n\} . \qquad (3.12)$$

Taking into account the stationarity of the increments of $\xi(t)$ and the following symmetry: $\xi(t) \doteq \xi(-t)$, one has

$$\xi_i - M_i \doteq M_\xi(D_i), \quad D_i = \{t_k - t_i, k = 0,1,...,i\} .$$

Suppose that $t_i = (t^{(1)},...,t^{(d)})_i \in S_n^2$, $t_i^{(1)} = -n$. Then, by (3.2),

$$D_i \supseteq [K_{\rho_n}], \quad K_{\rho_n} = [0, 2\rho_n] \times [-\rho_n, \rho_n]^{d-1},$$

where $\rho_n = n/L_n$.

Hence $M_\xi(D_i) \geq M_\xi([K_{\rho_n}])$ and

$$\max\{P(\xi_i - M_i \leq 0), t_i \in S_k^2, m \leq k \leq n\} \leq P(M_\xi([K_{\rho_{m,n}}]) \leq 0) , \qquad (3.13)$$

where

$$\rho_{m,n} = \min\{k/L_k, m \leq k \leq n\} \geq q\rho_n . \qquad (3.14)$$

The estimate (3.13) does not depend on the assumptions $t_i^{(1)} = -n$ since the process $\xi(t)$ is isotropic. It remains to combine (3.12), (3.13), and (2.7) to get

$$cn^{H-d}/\sqrt{\ln n} \leq P(M_\xi([K_{\rho_{m,n}}]) \leq 0) \leq P(M_\xi(q\rho_n \Delta) \leq \sqrt{\kappa \ln(q\rho_n)}), \qquad (3.15)$$

where $\Delta = [-1,1]^d + e = K_1(e)$ and $e = (1,0,...,0)$.

By setting

$$q\rho_n(\sqrt{\kappa \ln(q\rho_n)})^{1/H} = T, \quad \rho_n = n/(\ln n)^\varepsilon, \quad \varepsilon > 0$$

and assuming the self-similarity of $\xi(t)$, we get the desired relation for $\xi(t) = w_H(t)$:

$$P(M_{w_H}(T\Delta) \leq 1) \geq cT^{H-d}/(\ln T)^\lambda, \quad \lambda = d/(2H) + \varepsilon(d-H).$$

In the general case, (3.15) implies

$$cn^{H-d}/\sqrt{\ln n} \leq P(M_\xi([q\rho_n \Delta]) \leq 0)$$

or

$$P(M_\xi([T\Delta]) \leq 1) \geq cT^{H-d}/(\ln T)^\lambda, \quad \lambda = 1/2 + \varepsilon(d-H),$$

where $[T\Delta] = T\Delta \cap Z^d$.

**The upper bound**

To estimate $p_T(\Delta) = P\{\xi(t) < 1, t \in T\Delta\}$ from above it is sufficient to consider the spherical domain $\Delta = B_1(e) \subset R^d$ with center $e = (1,0,...,0)$. The desired estimate

$$p_T(\Delta) \leq \exp(c\sqrt{\ln T})T^{-d+H} \qquad (3.16)$$

for the fractional Brownian motion was obtained in the paper ( Molchan-1, 2017). The proof uses only two facts:

(a) $EM_\xi(T\Delta) \sim cT^H EM_\xi(\Delta)$ and

(b) existence of an element $\varphi(t) \in H_\xi(R^d)$ such that $\varphi(t) = 1, |t| \geq 1$.

The sufficient conditions for (a,b) are contained in Lemmas 4-5. Therefore, the upper bound (3.16) can be considered as proved.

.

**Acknowledgements.** This research was supported by the Russian Science Foundation through the research project 17-11-01052.